\documentstyle[12pt]{article}
\newcommand{\fr}{ \frac}

\newcommand{\Gc}{ {\cal G}}
\newcommand{\Jc}{ {\cal J}}
\newcommand{\Lc}{ {\cal L}}
\newcommand{\Zb}{\mbox{\boldmath $Z$ } }

\newcommand{\beq}{\begin{equation}}
\newcommand{\eeq}{\end{equation}}
\newcommand{\lb}{\label}

\newcommand{\ug}{\upsilon}

\newcommand{\ar}{\rightarrow}

\begin{document}

\vspace{2cm}
\begin{flushright}
FGI-13-98 \\
math.QA/9812108
\end{flushright}

\begin{center}
{\Large {\bf Green Function on the Quantum Plane}}
\end{center}

\vspace{1cm}

\begin{flushleft}
H. Ahmedov$^1$  and I. H. Duru$^{2,1}$
\vspace{.5cm}

{\small 
1. Feza G\"ursey Institute,  P.O. Box 6, 81220,  \c{C}engelk\"{o}y, 
Istanbul, Turkey 
\footnote{E--mail : hagi@gursey.gov.tr and duru@gursey.gov.tr}.

2. Trakya University, Mathematics Department, P.O. Box 126, 
Edirne, Turkey.}
\end{flushleft}

\vspace{1cm} \noindent 
{\bf Abstract}: Green function (which can be called the q-analogous of the
Hankel function) on the quantum plane $E_q^2=  E_q(2)/U(1)$ is constructed.

\section{Introduction}

Green functions play important roles in physics. Field theoretical problems
involving boundaries, such as the Casimir interactions, particle pair
productions i.e., all employ Green functions. Therefore if one is interested
in the investigation of some physical effects on the non--commutative spaces
construction of the Green functions in these media is useful. Motivated by
these considerations we think it is of interest to study the Green functions
on the quantum group spaces which are the natural examples of the
non--commutative geometries.

Previously we have constructed the Green function on the quantum sphere $%
S_q^2$ \cite{1}. In this paper we study the same problem for the quantum
plane $E_q^2$ which may be more relevant to physics.

In Section 2 we recall main result \cite{2,3,4,5,6,7} concerning the quantum
group $ E_q(2)$ and its homogeneous spaces.

In Section 3 we construct the Green function on the  quantum plane. The
Green function we obtain, provides the possibility of the future studies on
the new q-functions which are the deformations of the Neumann and Hankel
functions.

\newpage

\section{Quantum Group $E_q(2)$ and its Homogeneous Spaces}

Let $A$ be the set of linear operators in the Hilbert space $l^2(%
\mbox{\boldmath $Z$ } )$ subject to the condition 
\begin{equation}
\label{1}\sum_{j=-\infty}^\infty q^{2j} (e_j, F^*F e_j)< \infty ; \ \ \ \ \
\ \ F\in A. 
\end{equation}
Here $0< \ q < 1$ and $\{e_j \}$ is the orthonormal basis in $l^2(%
\mbox{\boldmath $Z$ } )$. Explicit form of $e_j$ is 
\begin{equation}
e_j=(0, \cdots , 0,\ 1, \ 0, \cdots ), 
\end{equation}
where either $j^{th}$ ( for $j>0$ ) or $(-\mid j\mid )^{th}$ ( for $j<0$ )
component is one, all others are zero. Any vector $x=(x_0, \ x_1, \ x_{-1},
\cdots , x_n, \ x_{-n}, \cdots )$ of $l^2(\mbox{\boldmath $Z$ } )$ has
representation 
\begin{equation}
x=\sum_{j=-\infty}^\infty x_je_j. 
\end{equation}
$(\cdot , \cdot )$ in (\ref{1}) is the scalar product in $l^2(%
\mbox{\boldmath $Z$ } )$: 
\begin{equation}
(x,y)=\sum_{j=-\infty}^\infty \overline{x_j}y_j. 
\end{equation}
$A$ is the Hilbert space with the scalar product 
\begin{equation}
\label{sca}(F, G)_A= (1-q^2)\sum_{j=-\infty}^\infty q^{2j} (e_j, F^*G e_j);
\ \ \ \ \ \ \ F, \ G\in A. 
\end{equation}

\vspace{3mm} \noindent
Let us introduce the linear operators acting in $l^2(\mbox{\boldmath $Z$ } )$
\begin{equation}
ze_j =e^{i\psi} q^j e_j, \ \ \ \ \upsilon e_j =e^{i\phi} e_{j+1}, 
\end{equation}
where $\psi$ and $\phi$ are the classical phase variables. $n$ is normal and 
$\upsilon$ is unitary operator in $l^2(\mbox{\boldmath $Z$ } )$. It is easy
to show that they satisfy the relations : 
\begin{equation}
z\upsilon =q\upsilon z, \ \ \ z^*\upsilon =q\upsilon n^*, \ \ \ zz^* =z^*z. 
\end{equation}
Any element $F\in A$ can be represented as 
\begin{equation}
F=\sum_{j=-\infty}^\infty f_j(z,z^*)\upsilon^j 
\end{equation}
by suitable choice of the functions $f_j$.

\vspace{3mm} \noindent
The linear operators $Z$ and $V$ given by 
\begin{equation}
Z=z\otimes\upsilon^{-1}+\upsilon\otimes z, \ \ \ V=\upsilon\otimes \upsilon 
\end{equation}
are normal and unitary in $l^2(\mbox{\boldmath $Z$ } \times 
\mbox{\boldmath $Z$ } )$ They satisfy the relations 
\begin{equation}
ZV =qVZ, \ \ \ Z^*V =qV Z^*, \ \ \ ZZ^* =Z^*Z. 
\end{equation}
Note that the operators $N$ and $V$ have the same properties as $n$ and $%
\upsilon$. Therefore there exits the linear map \cite{2} 
\begin{equation}
\Delta : A\rightarrow A\otimes_A A 
\end{equation}
defined as 
\begin{equation}
\Delta (f(z, z^*)\upsilon^j) = f(Z, Z^*)V^j. 
\end{equation}
Here $\otimes_A$ is the completed tensor product $\otimes$ with respect to
the scalar product 
\begin{equation}
(F_1\otimes F_2, F_3\otimes F_4)_A = (F_1, F_3)_A(F_2, F_4)_A; \ \ \ \
F_n\in A 
\end{equation}
in $A\otimes A$. $A$ is the space of square integrable functions on the
quantum group $E_q(2)$ and $\Delta$ is the quantum analog of the group
multiplication.

\vspace{3mm} \noindent
The one parameter groups $\{\sigma _1\}$ and $\{\sigma _2\}$ of automorphism
of $A$ given by 
\begin{equation}
\sigma _1(\upsilon )=e^{-it}\upsilon ,\ \ \ \sigma _1(z)=e^{it}z
\end{equation}
and 
\begin{equation}
\sigma _2(\upsilon )=\upsilon ,\ \ \ \sigma _2(z)=e^{it}z
\end{equation}
with $t\in \mbox{\boldmath $R$}$ are isomorphic to $U(1)$. The subspaces 
\begin{equation}
B=\{F\in A:\sigma _1(F)=F,\ \ {\rm for\ all}\ \ t\in \mbox{\boldmath $R$}\}
\end{equation}
and 
\begin{equation}
H=\{F\in B:\sigma _2(F)=F\ \ {\rm for\ all}\ \ t\in \mbox{\boldmath $R$}\}
\end{equation}
are the space of square integrable functions on the quantum plane $E_q^2$
and two sided coset space $U(1)\backslash E_q(2)/U(1)$. Any element of $H$
is the function of $\rho =zz^{*}$. Note that the scalar product (\ref{sca})
on $H$ becomes a q-integration 
\begin{equation}
\label{sca1}(f(\rho ),g(\rho ))_A=(1-q^2)\sum_{j=-\infty }^\infty q^{2j}%
\overline{f(q^{2j})}g(q^{2j})=\int_0^\infty \overline{f(\rho )}g(\rho
)d_{q^2}\rho 
\end{equation}
\vspace{3mm} \noindent 
Let $U_q(e(2))$ be the $*$--algebra generated by $p$ and $\kappa ^{\pm 1}$
such that 
\begin{equation}
p^{*}p=q^2pp^{*},\ \ \ \kappa ^{*}=\kappa ,\ \ \ \kappa p=q^2p\kappa .
\end{equation}
The representation ${\cal L}$ of $U_q(e(2))$ in $A$ is given by 
\begin{eqnarray}\lb{rep}
\Lc (p) f(z,z^*)\ug^j  & = & iq^{j+1}  D^z_+  f(z,z^*)\ug^{j+1}  \\
\Lc (p^*) f(z,z^*)\ug^j   & = & iq^j D^{z^*}_- f(z,z^*)\ug^{j-1} \\
\Lc (\kappa ) f(z,z^*)\ug^j   & = & q^j  f(q^{-1}z, qz^*)\ug^j, 
\end{eqnarray}
where 
\begin{equation}
D_{\pm }^xf(x)=\frac{f(x)-f(q^{\pm 2}x)}{(1-q^{\pm 2})x}.
\end{equation}
For the Casimir element $C=-q^{-1}\kappa^{-1}pp^{*}$ we have 
\begin{equation}
\Lc (C)f(z,z^{*})\upsilon^j=q^j D_{-}^{z^{*}}D_{+}^zf(qz, q^{-1}z^{*})
\upsilon ^j.
\end{equation}
The restriction $\Box $ of ${\cal L}(C)$ on $H$ is 
\begin{equation}
\Box =D_{-}^\rho \rho D_{+}^\rho 
\end{equation}
which we call the radial part of ${\cal L}(C)$.

\newpage

\section{Green Function on the Quantum Plane}
{\bf (i) Green Function on $U(1)\backslash E_q(2)/U(1)$}
\vspace{2mm}
\noindent

The Green function ${\cal G}^p(\rho )$ on the two sided coset space is
defined as 
\begin{equation}
\label{eqs}(\Box + p){\cal G}^p(\rho )=\delta (\rho ),
\end{equation}
where $\delta $ is the delta function which defined with respect to the
scalar product (\ref{sca1}) as 
\begin{equation}
(\delta ,f)_A=f(0)
\end{equation}
for any $f\in H_0$. The equation (\ref{eqs}) is understood as 
\begin{equation}
({\cal G}_p,f)_A=\lim _{\epsilon \rightarrow 0}(\delta (\rho ),\frac 1{\Box
+p+i\epsilon }f(\rho ))_A
\end{equation}
For $\rho \neq 0$ the equation (\ref{eqs}) is solved by 
\begin{equation}
{\cal J}(\sqrt{p\rho})=\sum_{k=0}^\infty \frac{(-1)^k}{([k]!)^2}(p\rho )^k
\end{equation}
and 
\beq
{\cal N}(\sqrt{p\rho})=\frac{q-q^{-1}}{2q\log (q)}{\cal J}(\sqrt{p\rho})
(\log (p\rho)+2C_q) - \fr{1}{q}
\sum_{k=1}^\infty \frac{(-1)^k}{([k]!)^2}(p\rho )^k
\sum_{m=1}^k\frac{q^m+q^{-m}}{[m]},
\eeq
where 
\beq
[m]=\frac{q^m-q^{-m}}{q-q^{-1}},\ \ \ [m]!=[1][2]\cdots [m].
\eeq
The Hahn-Exton q-Bessel function ${\cal J}$ is regular at $\rho =0$. 
It is the zonal spherical function of the unitary irreducible 
representations of $E_q(2)$. 
${\cal N}$ can be called q-Neuman function which indeed is reduced to the
usual Neuman function in $q\rightarrow 1$ limit. 
Here $C_q$ is some constant, which in $q\rightarrow 1$ limit 
should become the Euler constant \cite{8}. 

The Green function on $U(1)\backslash E_q(2)/U(1)$ is then 
\begin{equation}
{\cal G}_p(\rho )={\cal N}(\sqrt{p\rho} )-i{\cal J}(\sqrt{p\rho}),
\end{equation}
which in classical limit becomes the Hankel function.
Using the Fourier-Bessel integral \cite{9}
\beq
\int_0^\infty d_{q^2}\rho \Jc (q^n\sqrt{\rho}) \Jc (q^m\sqrt{\rho}) = 
\fr{q^{2m+2}}{1-q^2} \delta_{mn}
\eeq
we arrive at the following representation for the Green function 
\beq
\Gc_p (\rho ) = \lim _{\epsilon \rightarrow 0}
q^{-2}\int_0^\infty d_{q^2}\lambda \fr{\Jc (\sqrt{\lambda\rho} )}
{p-\lambda+i\epsilon }
\eeq
from which one can derive the constant $C_q$.

\vspace{2mm} \noindent
To prove that $ {\cal G}$ solves (\ref{eqs}) we first have to show that 
\begin{equation}
\Box \log \rho = \frac{2q\log (q) }{q-q^{-1}}\delta (\rho ). 
\end{equation}
For $\rho\neq 0$ we have 
\begin{equation}
\Box \log \rho =0. 
\end{equation}
Since the operator $\Box$ is symmetric in $H$ we have 
\begin{eqnarray}
( \Box \log \rho , \ f )_A & = & 
( \log \rho , \ \Box f )_A \nonumber \\
& = & \frac{2q\log (q) }{q-q^{-1}}\lim_{n\ar\infty}\sum_{j=-\infty}^n 
j [ 2f(q^{2j}) - f(q^{2(j+1)})-f(q^{2(j-1)})] \nonumber \\
& = &\frac{2q\log (q) }{q-q^{-1}}\lim_{n\ar\infty}[f(q^{2n}) -
n( f(q^{2n+2})- f(q^{2n})) ].
\end{eqnarray}
We then employ the q--Taylor expansion at the neighborhood of $\rho =0$ 
\begin{equation}
f(q^2\rho )-f(\rho )\sim D_+^\rho f(0) (q^2- 1)\rho . 
\end{equation}
For $n>> 1$ we get 
\begin{equation}
n( f(q^{2n+2})- f(q^{2n}))\sim n D_+^\rho f(0) (q^2- 1)q^{2n}. 
\end{equation}
Since $nq^{2n} $ vanishes as $n\rightarrow\infty$, we arrive at 
\begin{equation}
(\Box \log \rho , \ f )_A = \frac{2q\log (q) }{q-q^{-1}}\lim_{n\rightarrow%
\infty}f(q^{2n}) = \frac{2q\log (q) }{q-q^{-1}}f(0). 
\end{equation}
In a similar fashion one can show that 
\begin{equation}
( (\Box+p)  {\cal G}^p (\rho ), \ f )_A = f(0). 
\end{equation}

\vspace{4mm}
\noindent
{\bf (ii) Green Function on $E_q^2$}
\vspace{2mm}
\noindent

We obtain the Green function $ {\cal G}^p (R )$ on the quantum plane $E_q^2$
from the one $ {\cal G}^p (\rho )$ on the two sided coset space by the group
multiplication \cite{1} : 
\begin{equation}
{\cal G}^p (R )= \Delta  {\cal G}^p (\rho ). 
\end{equation}
Here 
\begin{equation}
R= \Delta (\rho )= \rho \otimes 1 + 1\otimes \rho + \upsilon z^*\otimes
z\upsilon + z\upsilon^* \otimes \upsilon^* z^* 
\end{equation}
is self-adjoint operator in $l^2(\mbox{\boldmath $Z$ } \times 
\mbox{\boldmath $Z$ } )$ and 
\begin{equation}
R e_{ts}=q^{2t}e_{ts} 
\end{equation}
where the eigenfunctions $e_{ts}$ are given by \cite{3} 
\begin{equation}
e_{ts}= \sum_{j=-\infty}^\infty (-1)^j q^{t-j} {\cal J}_s
(q^{t-j})e_{s+j}\otimes e_j . 
\end{equation}
They satisfy the orthogonality condition 
\begin{equation}
(e_{ts},e_{ij})=\delta_{ti}\delta_{sj}. 
\end{equation}
We also have 
\begin{equation}
e_{s+j}\otimes e_j =\sum_{t=-\infty}^\infty (-1)^jq^{t-j} {\cal J}_s
(q^{t-j})e_{ts}. 
\end{equation}
Therefore the basis elements $e_{ts}$; $t,s\in (-\infty, \ \infty )$ form
the complete set in $l^2(\Zb  \times \Zb )$. 
The Green function on the quantum plane is the linear operator in this
space defined as 
\begin{equation}
{\cal G}^p (R )e_{ts}=  {\cal G}^p (q^{2t})e_{ts}. 
\end{equation}


\end{document}